\documentclass[10pt]{article}
\usepackage[english]{babel}
\usepackage{graphicx}
\usepackage{framed}
\usepackage[normalem]{ulem}
\usepackage{amsmath}
\usepackage{amsthm}
\usepackage{amssymb}
\usepackage{amsfonts}
\usepackage{enumerate}
\usepackage[utf8]{inputenc}
\usepackage[top=1 in,bottom=1in, left=1.75 in, right=1.75 in]{geometry}

\title{A translation of ``On the positive representation of polynomials" by Ernst Meissner}
\author{Hannah K. Wayment-Steele}
\date{last updated 05/27/2023}
\begin{document}

\maketitle

\section*{Translator's note}
We provide here an English translation of

\hfill

\noindent Meissner, Ernst. \"Uber positive Darstellungen von Polynomen (1911) \emph{Mathematische Annalen} 70, 223-235.

\hfill

\noindent We used Meissner's result that any polynomial may be written as a rational polynomial (Equation 3) to motivate our exploration of biomolecules and their partition functions as universal function approximators.

\hfill

\noindent We thank Rhiju Das and Cole Graham for useful discussion. We welcome feedback on the translation to wayment (at) brandeis (dot) edu. We assert that the original work is in the public domain via the German Act on Copyright and Related Rights, Section 64 (copyright expires 70 years after the author’s death). 

\section{On the Positive Representation of Polynomials}

In the present article I will be occupied with the following theorem:

\hfill

Theorem: Let
\begin{equation}
 f(x) = x^k + c_1 x^{k-1} + \cdots + c_k   
\end{equation}
be a polynomial with real coefficients that satisfies the condition
\begin{equation}
    f(x) > 0, x \geq 0.
\end{equation}

Then $f(x)$ can always be written in the form
\begin{equation}
    f(x) = \dfrac{f_2(x)}{f_1(x)}
\end{equation}

where $f_1(x)$ and $f_2(x)$ are positive-coefficient polynomials.

\hfill

The representation (3) will be called a \emph{positive representation of polynomial $f(x)$}. The properties stated in (2) are critical.

This theorem is said to result from Laguerre, but it is unknown to me where and how Laguerre demonstrated it.\footnote{It does not appear that the works of Laguerre contain the theorem. [All footnotes are from EM.]} Herr A.~Hurwitz presented in a seminar in 1904 how one may achieve a positive representation through a special approach (see footnote [2]).

In the following, I will develop a very general method to construct a positive representation.  It follows a reduction of the problem that I adopted from Mr.~Hurwitz, essentially the solution of a linear system of inequalities with positive variables. My method depends on the appropriate geometric interpretation of the problem via the construction of a convex polygon, and is carried out with a convex curve.

The second part of my work will address the question of the extension of the theorem to polynomials of two variables (and therefore in principle any number). It is noteworthy that a positive representation can only be reached if one tightens the constraints in (2).

\subsection{The Hurwitzian reduction of the problem to the solution of a system of inequalities.}

Let $f(x)$ be a polynomial in (1) with real coefficients, that satisfies the requirements in (2); one can separate it into real and real irreducible (linear and quadratic) factors.  Then one only needs to prove the Laguerrian theorem for these, as the product of positive representations of two polynomials provides a positive representation of their product.

Because of (2), the linear real factors of $f(x)$ have the negative roots
\begin{equation*}
    x_i = -|x_i|,
\end{equation*}
which in the equation
\begin{equation*}
    f(x)=0
\end{equation*}
take the form
\begin{equation*}
    (x-x_i) = (x+|x_i|).
\end{equation*}
They therefore appear in the positive representation without further modification.

The quadratic factors have complex conjugate roots, which we write as
\begin{equation*}
    x=re^{\pm i \phi}, r > 0, 0 < \phi < \pi.
\end{equation*}

These are therefore of the form
\begin{equation*}
    (x-re^{i\phi})(x-re^{-i\phi}) = x^2 - 2 x r \cos \phi + r^2.
\end{equation*}
All the coefficients are also positive as long as
\begin{equation*}
    \pi/2 < \phi < \pi.
\end{equation*}
Therefore in the following we restrict ourselves to the case
\begin{equation*}
    0 < \phi < \pi/2.
\end{equation*}
To exact a positive representation for the factor
\begin{equation*}
    Q(x) = x^2 - 2 x r \cos \phi + r^2,
\end{equation*}
we must find a polynomial with only positive coefficients $Q_1(x)$, such that the polynomial $Q_2(x) = Q(x)Q_1(x)$ also has positive coefficients. We set $Q_1(x)$ and $Q_2(x)$ as

\begin{equation}
\begin{split}
Q_1(x) &= a_0 + \dfrac{a_1}{r}x + \dfrac{a_2}{r^2}x^2 + \cdots + \dfrac{a_n}{r^n}x^n\\   
Q_2(x) &= r^2\left(a_0 + \dfrac{b_0}{r}x + \dfrac{b_1}{r^2}x^2 + \cdots + \dfrac{b_n}{r^{n+1}}x^{n+1} + \dfrac{a_n}{r^{n+2}}x^{n+2}\right).\\ 
\end{split}
\end{equation}

Here the constants $a$ are positive and so chosen that the constants $b$ are also positive. We set
 \begin{equation*}
a_{-1} = 0, \, a_{n+1} = 0,
 \end{equation*}
 so we have in general, for $k=0,1,\ldots n$,
  \begin{equation*}
  b_k = a_{k-1} + a_{k+1} - 2a_k \cos \phi
 \end{equation*}
 and we have therefore reduced the task to solving the following system of inequalities:
 \begin{equation}
     \begin{split}
         a_{k-1} + a_{k+1} - 2a_k \cos \phi &> 0,\\
         a_k &> 0,\\
         a_{-1} = 0, \, a_{n+1} &= 0.
     \end{split}
 \end{equation}
 
 \subsection{Geometric Interpretation and Solution to the System of Inequalities.}
 In a right-handed coordinate system $O(x,y)$, $g(p,u)$ represents the line that lies a distance $p$ from the origin $O$ and makes the angle $u$ with the x-axis.
 
 Let $a_{k-1}, a_k, a_{k+1}$ be three positive constants, and $\phi$ an angle between $0$ and $\pi/2$.  We draw the three lines
 \begin{equation*}
 \begin{split}
     g_{k-1} &= g(a_{k-1}, \, k\phi),\\
     g_k &= g(a_k, \,(k+1)\phi),\\
     g_{k+1} &= g(a_{k+1},\,(k+2)\phi).
 \end{split}
 \end{equation*}
 The line segment of $g_k$ that is defined by intersecting with lines $g_{k-1}, g_{k+1}$ has length
  \begin{equation*}
\dfrac{1}{\sin \phi}(a_{k-1}+a_{k+1}-2a_k\cos \phi) = \dfrac{b_k}{\sin \phi}.
 \end{equation*}
 
This formula expresses this length as either positive or negative, depending on if the origin $O$ lies on the concave or convex side of the line drawn by the three lines (Fig.~1 or Fig.~2, respectively).

\includegraphics[width=0.5\textwidth]{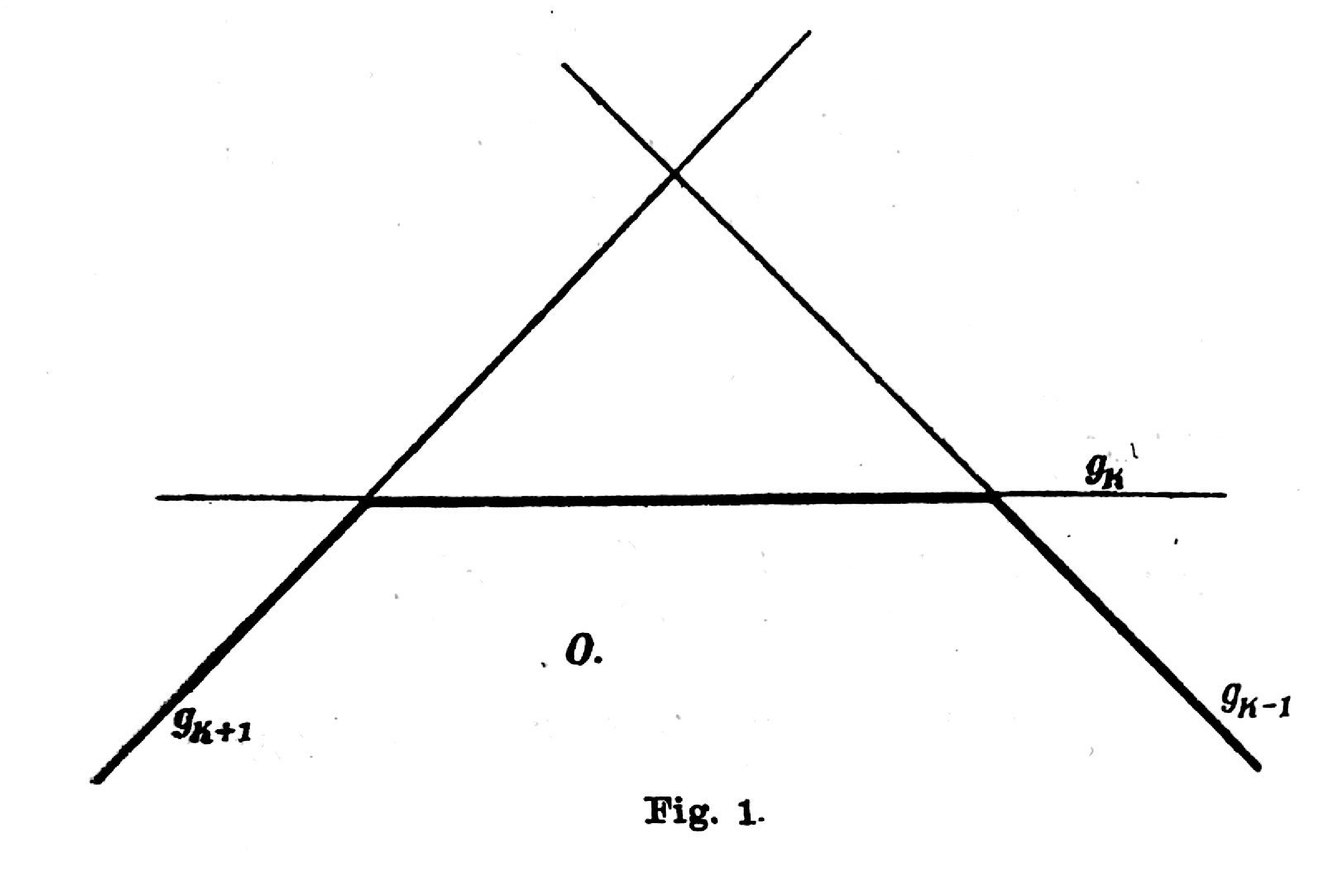}
\includegraphics[width=0.5\textwidth]{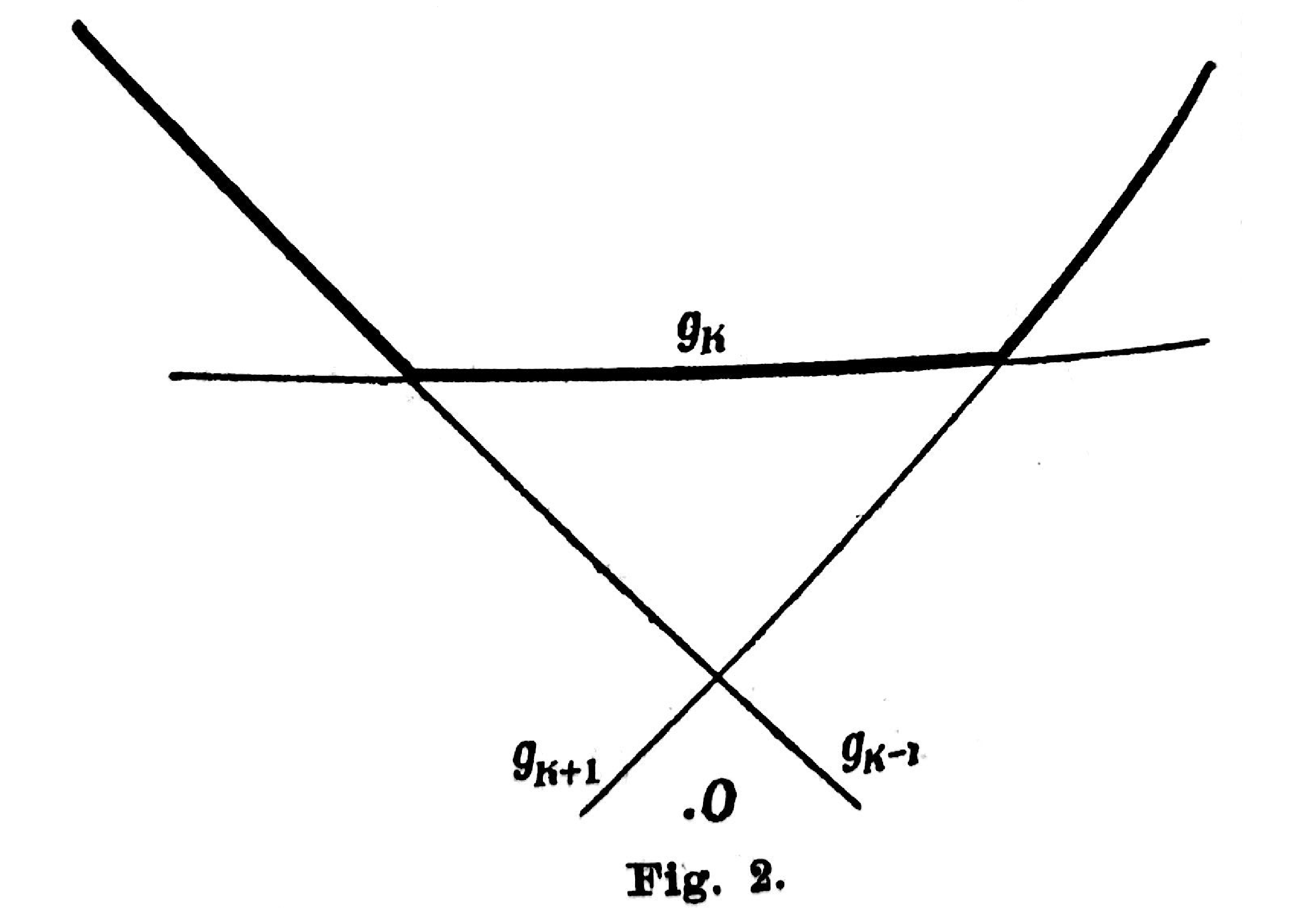}
 
 If the expression is $\dfrac{1}{\sin \phi}(a_{k-1}+a_{k+1}-2a_k\cos \phi) = 0$, then the three lines go through one point.
 
Consider a system of constants $a_k, \phi$ that satisfy the inequalities in (5), and construct the according lines
  \begin{equation*}
     g_k = g(a_k, (k+1)\phi), k = -1, 0, 1, \ldots, n+1.
 \end{equation*}
Two successive straight lines of the polygon form the angle $\phi$, with the first and the last going through the origin $O$, since $a_{-1} = a_{n+1} = 0$. The inequalities
  \begin{equation*}
a_{k-1}+a_{k+1}-2a_k\cos \phi > 0
 \end{equation*}
 now simply express the fact that this equiangular polygon $P(\phi)$ is concave toward the origin everywhere.
 It is evident that every such equiangular polygon that is concave against the origin, starts and ends at the origin, and has the distances $a_k$ from the origin is a solution to the inequality system in (5).
 
 All solutions of (5) result from the construction of all polygons $P(\phi)$. This leaves no further difficulties.\footnote{This immediately results in the minimal value of $n$ for which (5) has solutions for a given angle $\phi$.  One must draw the polygon $P(\phi)$ so that it is possible to return to the origin. This depends on how small the sides of $P(\phi)$ may possibly be. If $b_k$ may equal $0$, then all the lines $g_1, g_2, \ldots, g_n$ go through the intersection $S$ of $g_{-1}$ and $g_0$. The minimal value of $n$ then is written as
  \begin{equation*}
(n+1)\phi < \pi \leq (n+2)\phi;
\end{equation*}
as then the line $g_{n+1}$ makes the angle $(n+2)\phi > \pi$ with the x-axis, and the first line from $P(\phi)$, which may go through $S$ and the origin, without allowing $b_k$ to be negative (Fig.~3). 

\begin{center}
 \includegraphics[width=0.75\textwidth]{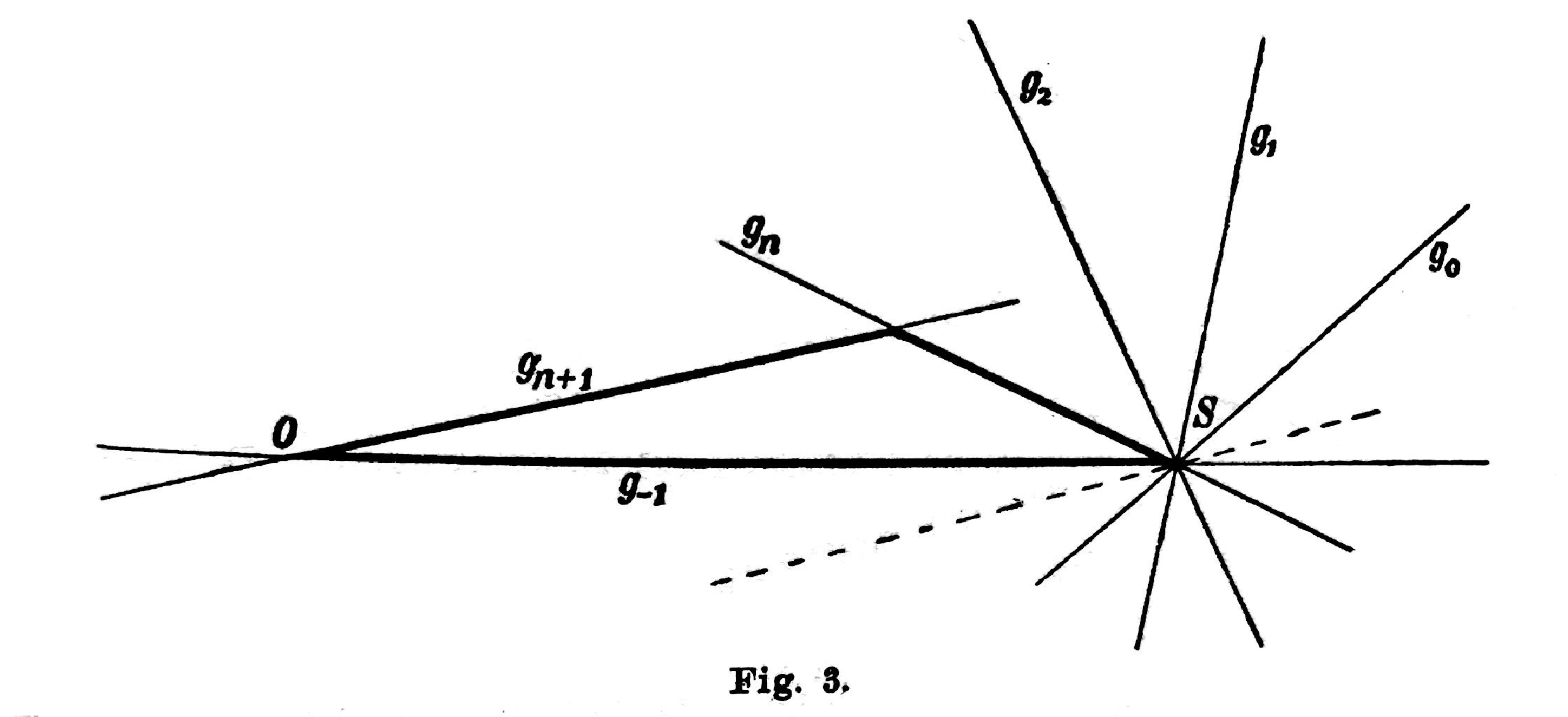}
\end{center}
The polynomial $P(\phi)$ is reduced to the triangle $(g_{-1}, g_n, g_{n+1})$, with one corner at $S$ which is counted multiple times. Using an arbitrary choice for $a_0$, the coefficients $a_k$ are written as
  \begin{equation*}
a_k = \dfrac{a_0}{\sin \phi} \sin\left((k+1)\phi\right), k = 1\ldots n.
\end{equation*}
 In the polynomial in (4), $Q_2(x)$ retains only the first and the last two coefficients as nonzero. This is the specific solution from Mr.~Hurwitz mentioned in the introduction.
 
 If the assumption $b_k = 0$ is invalid, then the solutions with minimal value $n$ may be obtained as follows. One only has to pull apart the corners of P that overlap in S, so that $b_k > 0$. The minimal value of $n$ does not change generally (Fig.~4).

 \begin{center}
  \includegraphics[width=0.75\textwidth]{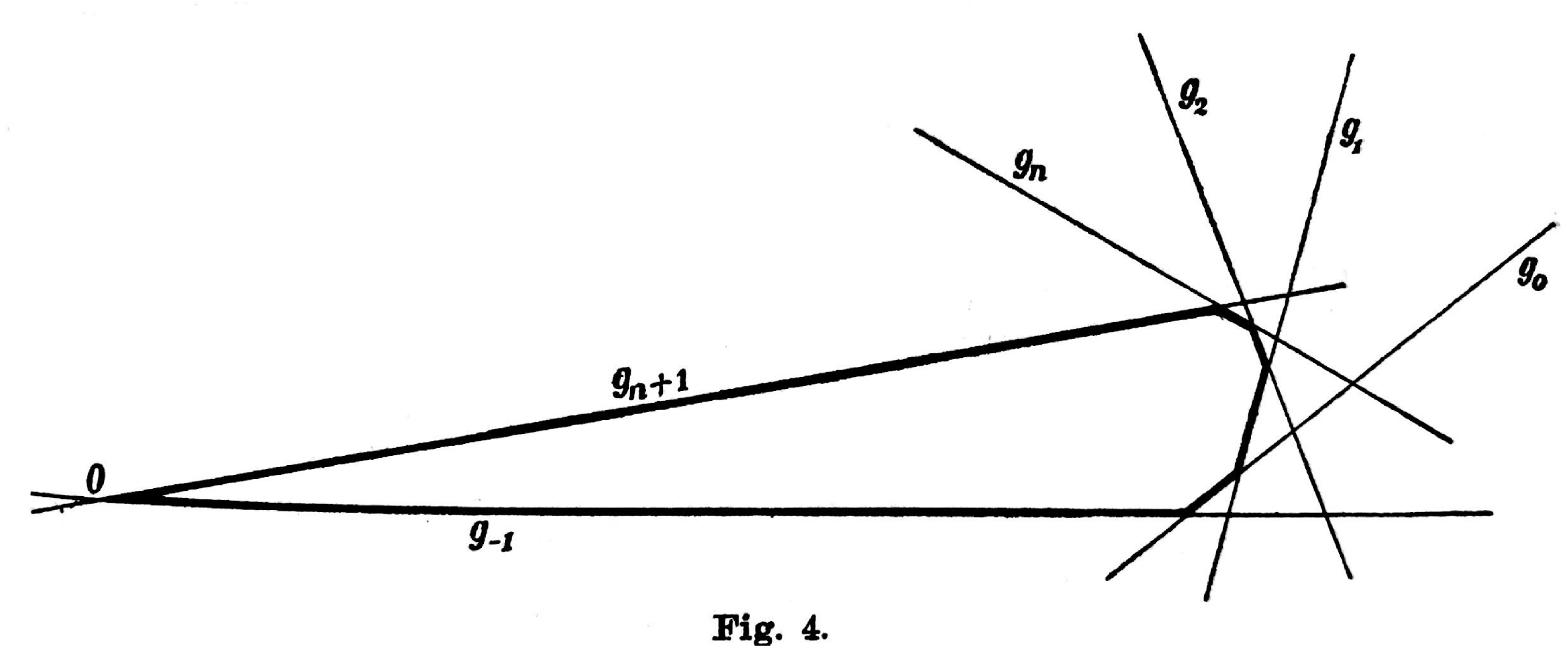}
\end{center}
 Only if $\pi/\phi$ is an integer (in which case the above triangle turns into the line segment $OS$), does one need to increment $n$ in order to make all $b_k > 0$. The minimal value $n$, for which the system in (5) has solutions, may be therefore written via the inequalities
   \begin{equation*}
(n+1)\phi \leq \pi < (n+2)\phi.
\end{equation*}}

For every polygon $P(\phi)$, one can draw a curve $C$ inside it with the following properties:

a) $C$ touches all sides of $P(\phi)$ in turn,

b) $C$ has a radius of curvature that is a continuous function of the arc length, and does not become zero or infinity anywhere.

Conversely, consider a curve $C$ with the property b). One can then select the point $O$, such that at least two tangents $t_1$, $t_2$ go from $O$ to the curve $C$ in the manner shown in Fig.~5.

\begin{center}
\includegraphics[width=0.5\textwidth]{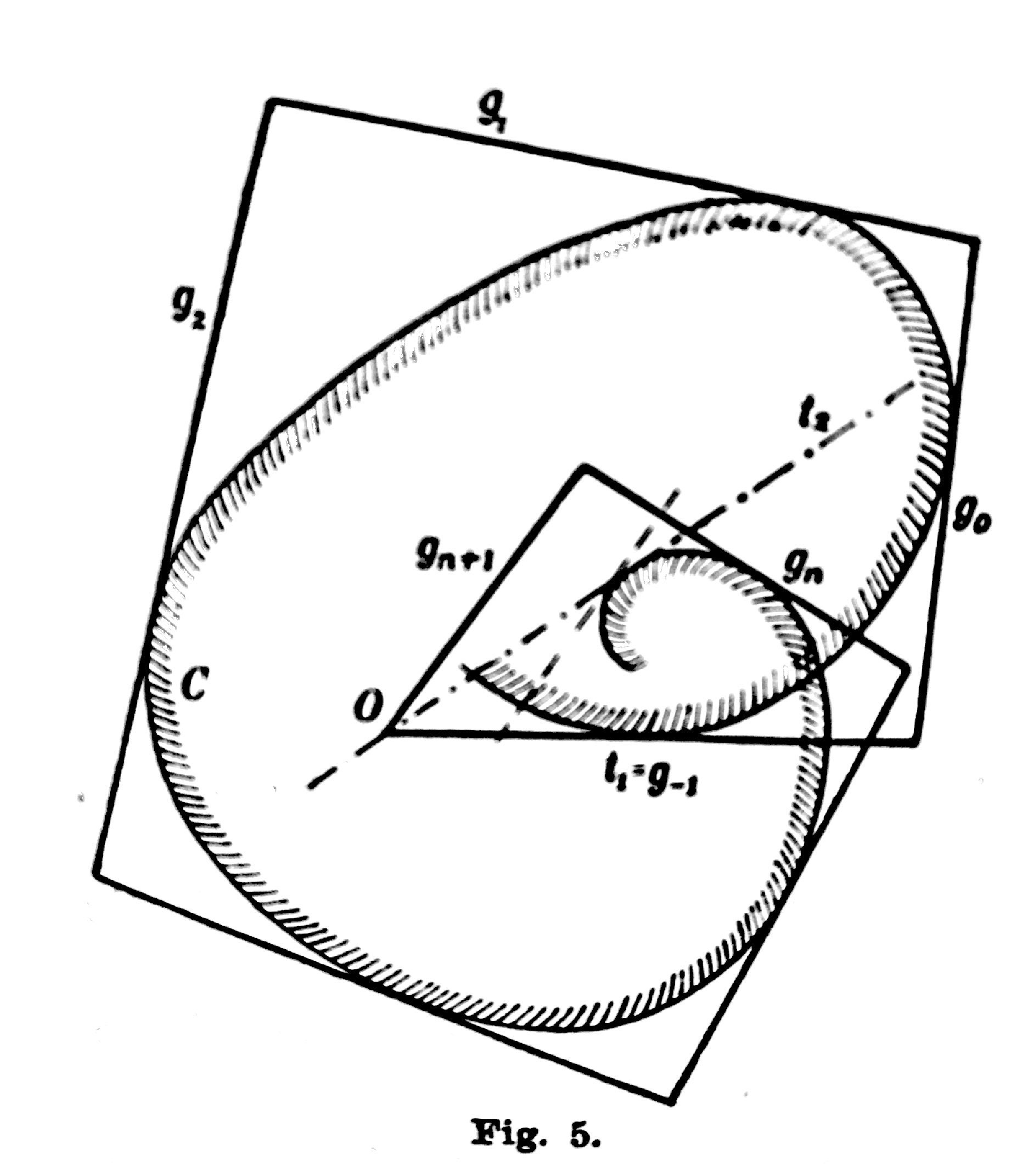}
\end{center}

Starting from $t_1$, one may construct an equiangular polygon $P(\phi)$ around $C$ with the polygon angle $\phi$. The first side of the polygon that extends beyond the tangent $t_2$ is replaced by a parallel line $g_{n+1}$ through O, and all subsequent sides are removed. Then one has a polygon which is a solution to the system in (5).  If the angle $\phi$ takes a second value $\phi_1$, then the same curve $C$ may be used in the same way to construct an associated solution to (5).  If the acute angle $\phi_1$ is greater than $\phi$, then it is clear that the solution according to $\phi$ also provides a solution for the value $\phi_1$.
\emph{All polygons $P(\phi)$ and therefore all solutions of the system (5) can be produced by drawing equi-angled polygons of all curves C in the specified manner.}

We will formulate this result analytically.

Let $t(u)$ be the tangent to $C$ that encloses the angle $u$ with $t_1$; $T(u)$ their contact point; and $p(u)$ their distance from point $O$. The function $p(u)$, the \emph{Support line function of curve $C$}, characterizes its form. The radius of curvature $\rho(u)$ of the curve $C$ at contact point $T(u)$ is given by the expression
\begin{equation}
    \rho(u) = p(u) + \dfrac{d^2 p(u)}{du^2}.
\end{equation}
For C to satisfy our requirements, we determine $p(u)$ from the differential equation (6) assuming $\rho(u)$ is an arbitrary continuous positive function.

We assert the two integration constants appearing in the integral so that $p(u)$ vanishes for a positive value $u_1$ and for $u = 0$, and is positive for intermediate values (Fig.~6). This is always possible.

\begin{center}
\includegraphics[width=0.8\textwidth]{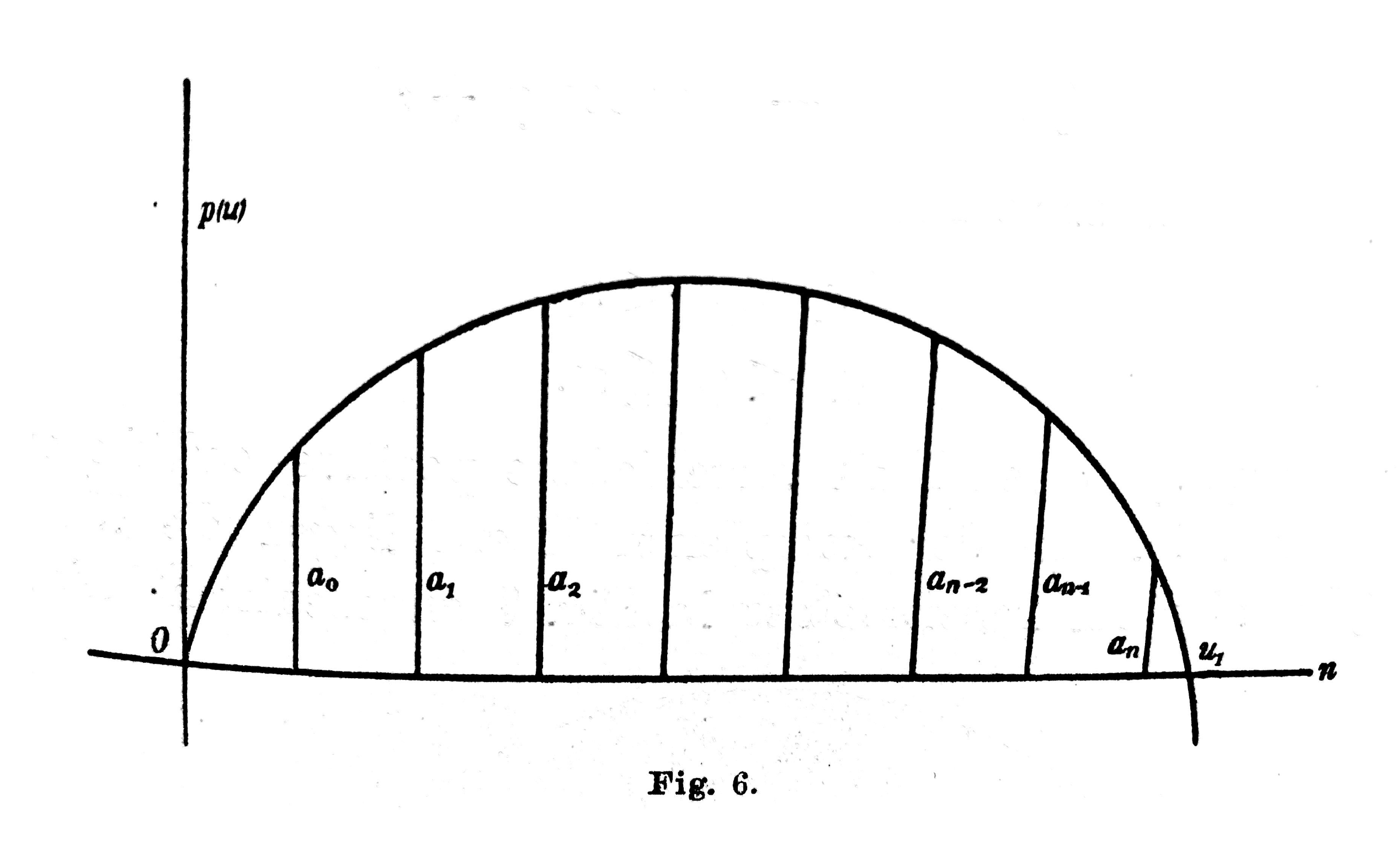}
\end{center}

\emph{The system of inequalities (5) is then solved through the formulas}
\begin{equation*}
    a_k = p((k+1)\phi), \, k = -1, 0, 1, \cdots, n
\end{equation*}
\emph{whereas the value of the integer $n$ comes from the relation}
\begin{equation*}
    (n+1)\phi < u_1 \leq (n+2) \phi.
\end{equation*}
This theorem demonstrates the perfect equivalence between the system of inequalities and the differential inequality\footnote{If the expression $a_{k-1} + a_{k+1} - 2a_k \cos \phi$ in (5) is allowed to be zero, then the system (5) is equivalent to the inequality 
\begin{equation*}
    p(u) + \dfrac{d^2 p(u)}{du^2} \geq 0.
\end{equation*}
The equality symbol here (and earlier) describes the solution with minimal $n$. The solution then becomes
\begin{equation*}
p(u) = \dfrac{a_0}{\sin \phi} \sin u, \text{ therefore } a_k = \dfrac{a_0}{\sin \phi} \sin (k+1) \phi.
\end{equation*}}
\begin{equation*}
    p(u) + \dfrac{d^2 p(u)}{du^2} > 0.
\end{equation*}

\subsection{Construction of the representation factors via a region $G$.}

To construct a positive representation of the original given polynomial $f(x)$, one must produce the two polynomials $Q_1(x), \, Q_2(x)$ with the specified method for each of its irreducible quadratic factors $Q(x)$ that do not have solely positive coefficients. For $Q(x)$ we thus have the positive representation
\begin{equation*}
    Q(x) = \dfrac{Q_2(x)}{Q_1(x)}.
\end{equation*}

These must be multiplied with themselves and with all other factors of $f(x)$. Thus $f(x)$ is expressed as a positive representation, via expression (3).

The polynomial $Q_1(x)$ will generally change from factor to factor as it depends on its roots. The work of this section is to remove this dependency.

Let the quantities $r$ and $\phi$ in the factor $Q(x) = x^2 - 2xr\cos\phi + r^2$ be subject to the constraints

\begin{equation}
    \begin{split}
        r_1 < r < r_2,\\
        \psi < \phi \leq \pi.\\
    \end{split}
\end{equation}
Here, $r_1,\, r_2, \, \psi$ are arbitrary nonnegative finite numbers. In the complex plane, the roots of the factor $Q(x)$ lie inside an open circular split region $G$ indicated in Fig.~7. For shorthand, we say $Q(x)$ belongs to region $G$.

\begin{center}
\includegraphics[width=0.8\textwidth]{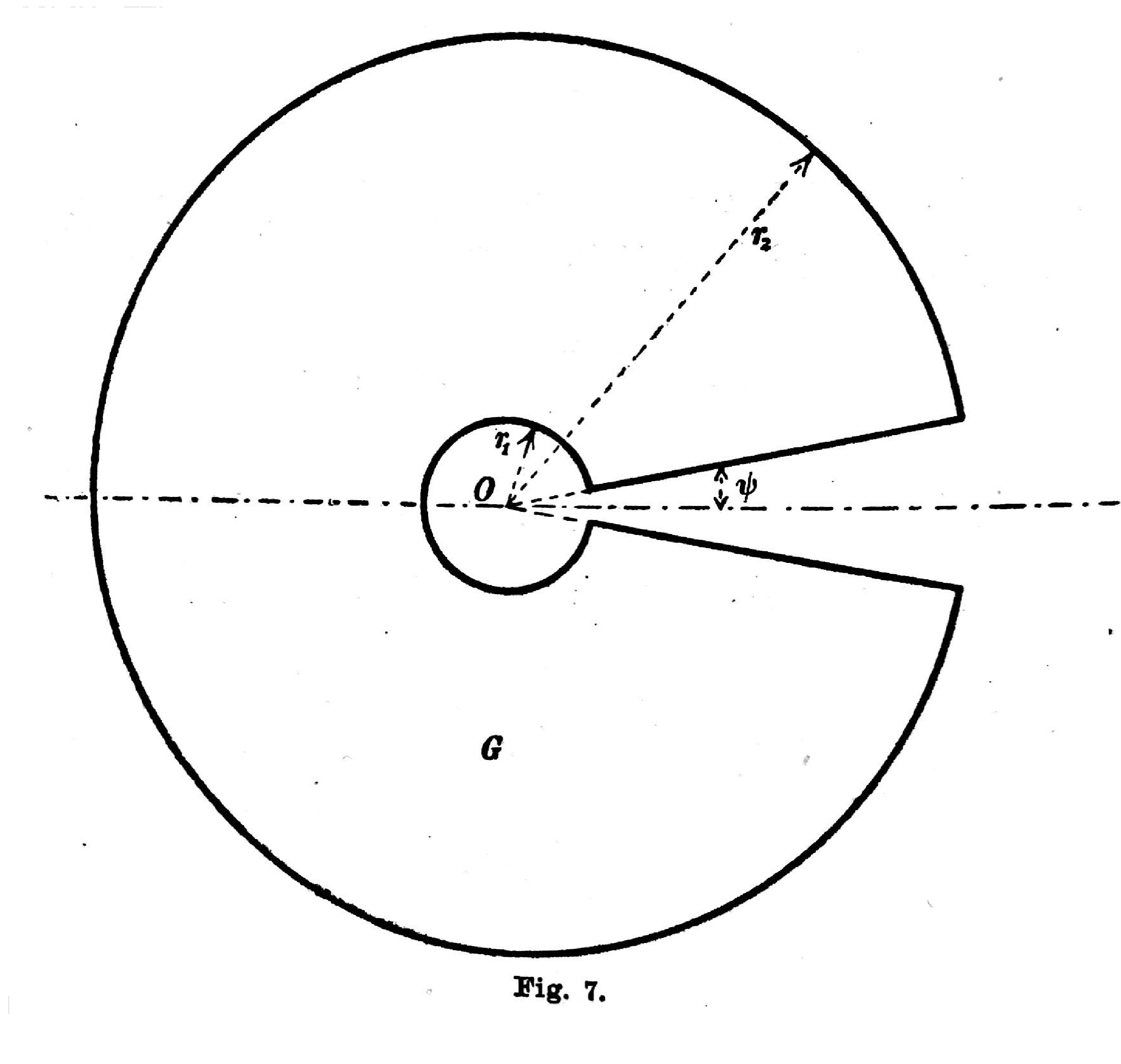}
\end{center}

\emph{Now we shall find a polynomial $G(x)$ with positive coefficients dependent only on the domain $G$, that is, on $r_1,\, r_2, \, \psi$, so that for all the factors $Q(x)$ belonging to the domain $G$, the polynomial}

\begin{equation*}
K(x) = Q(x) G(x)
\end{equation*}
\emph{will also have positive coefficients. We call $G(x)$ the representation factor of region $G$.}

We make the proposal
\begin{equation}
    \begin{split}
        G(x) &= \alpha_0 + \dfrac{\alpha_1}{r_2}x + \dfrac{\alpha_2}{r_2^2}x^2 + \cdots + \dfrac{\alpha_n}{r_2^n}x^n\\ 
        K(x) &= \alpha_0r^2 + \beta_0 rx + \beta_1 x^2 + \dfrac{\beta_2}{r_2}x^3 +\dfrac{\beta_3}{r_2^2}x^4 + \cdots + \dfrac{\beta_n}{r_2^{n-1}}x^{n+1} + \dfrac{\alpha_n}{r_2^{n}}x^{n+2}.
    \end{split}
\end{equation}
From here we use the abbreviation
\begin{equation*}
    \xi = \dfrac{r}{r_2}.
\end{equation*}
Then we have for the $\beta$ the equalities
\begin{equation*}
    \begin{split}
        \beta_0 &= \alpha_1 \xi - 2 \alpha_0 \cos \phi,\\
        \beta_1 &= \alpha_2 \xi^2 - 2 \alpha_1 \xi \cos \phi + \alpha_0,\\
        \beta_2 &= \alpha_3 \xi^2 - 2 \alpha_2 \xi \cos \phi + \alpha_1,\\
        &\ldots \\
        \beta_{n-1} &= \alpha_n \xi^2 - 2 \alpha_{n-1} \xi \cos \phi + \alpha_{n-2},\\
        \beta_n &=  - 2 \alpha_{n-1} \xi \cos \phi + \alpha_{n-1}.\\
    \end{split}
\end{equation*}
We set
\begin{equation*}
    \dfrac{r_1}{r_2} = \epsilon.
\end{equation*}
Due to (7) we have
\begin{equation*}
    0 < \epsilon < \xi < 1.
\end{equation*}
\emph{Now we have positive $\alpha_k$, and therefore need to ensure that all the inequalities}
\begin{equation*}
   \beta_k > 0
\end{equation*}
\emph{are fulfilled, and in fact for every value of $\xi$ in the interval} $\epsilon < \xi < 1$.
This is certainly the case if we set $\alpha_k$ positive and ensure the following inequalities:
\begin{equation}
    \begin{split}
        \gamma_0 &= \alpha_1 \xi - 2 \alpha_0 \cos \psi \geq 0,\\
        \gamma_1 &= \alpha_2 \xi^2 - 2 \alpha_1 \xi \cos \psi + \alpha_0 \geq 0,\\
        \gamma_2 &= \alpha_3 \xi^2 - 2 \alpha_2 \xi \cos \psi + \alpha_1 \geq 0,\\
        &\ldots \\
        \gamma_{n-1} &= \alpha_n \xi^2 - 2 \alpha_{n-1} \xi \cos \psi + \alpha_{n-2} \geq 0,\\
        \gamma_n &=  - 2 \alpha_{n-1} \xi \cos \psi + \alpha_{n-1} \geq 0,\\
        &(\epsilon < \xi < 1).\\
    \end{split}
\end{equation}
Then due to (7) the left-hand sides of these inequalities are smaller than their corresponding $\beta$ values.

We will not solve this system in general here, but content ourselves with specifying a particular solution that will reach for our purposes.

We arbitrarily select $\alpha_0 = 1$. Then $\alpha_1$ fits the inequality
\begin{equation*}
    \alpha_1 > \dfrac{2 \alpha_0 \cos \psi}{\epsilon} = \dfrac{2 \cos \psi}{\epsilon}.
\end{equation*}
Finally,
\begin{equation*}
 \alpha_k = \alpha_1^k \cos \psi^{k(k-1)}, \, k=2,3,\cdots n.
\end{equation*}

This statement holds for all integers $n$. Now all $\alpha$ are positive numbers. Furthermore,
\begin{equation*}
    \gamma_0 = \xi (\alpha_1 - \dfrac{2\cos\psi}{\xi}) > \xi (\alpha_1 - \dfrac{2 \cos \psi}{\epsilon}) > 0
\end{equation*}
and
\begin{equation*}
    \gamma_k = \alpha_{k-1}(\alpha_1 \xi \cos \psi^{2k-1} - 1)^2 \geq 0 \text{ for }k=1,2,\cdots n-1.
\end{equation*}
All coefficients $\gamma$ with the exception of the last one will therefore be positive, at least not negative. For the last one we have
\begin{equation*}
    \gamma_n = \alpha_{n-1}(1-2\alpha_1 \xi \cos \psi^{2n-1}).
\end{equation*}
If one chooses $n$ large enough so that
\begin{equation*}
    \cos \psi^{2n-1} < \dfrac{1}{2\alpha_1},
\end{equation*}
which is always possible because $\psi>0$, then $\gamma_n$ will also be positive and the system of inequalities in (9) is solved. Therefore the construction of the factors $G(x)$ is achieved.
One sees that the smaller $\psi$ is, the larger degree is chosen; also if the value of $\dfrac{r_1}{r_2} = \epsilon$ is small, the first coefficients $\alpha_1, \, \alpha_2, \ldots$ take large values. However, this is immaterial for the following.

\subsection{The Extension of the Laguerrian Theorem to Polynomials of Two Variables.}

We consider a positive representation of a polynomial $f(x,y)$ in two variables of the form
\begin{equation}
 f(x,y) = \dfrac{f_2(x,y)}{f_1(x,y)}
\end{equation}
where $f_1(x,y), \, f_2(x,y)$ are polynomials in $x$ and $y$ with nonnegative coefficients.

One might guess that, analogously to the Laguerrian theorem, such a representation exists, as long as 
\begin{equation}
    f(x,y) > 0 \text{ for } x>0,\, y>0.
\end{equation}

\emph{This is not the case.}

Let $f(x,y)$ have the positive representation in (10).  Let $\lambda$ and $\mu$ be the highest powers of $x$ in $f(x,y)$ and $f_1(x,y)$, respectively. If we take the limit of the following equation

\begin{equation*}
    \dfrac{f(x,y)}{x^\lambda} = \dfrac{\dfrac{f_2(x,y)}{x^{\lambda + \mu}}}{\dfrac{f_1(x,y)}{x^\mu}}
\end{equation*}

as $x = \infty$, and let $f(\infty, y)$ represent the following polynomial in $y$,

\begin{equation*}
L_{x=\infty} = \dfrac{f(x,y)}{x^\lambda},
\end{equation*}
in analogy to former expressions, thus follows 
\begin{equation*}
    f(\infty, y) = \dfrac{f_2(\infty, y)}{f_1(\infty, y)}.
\end{equation*}
The polynomials on the right hand side are nonzero with nonnegative coefficients. They therefore have positive values if $y$ is positive.
\emph{Every polynomial that has a positive representation also requires the constraint}
\begin{equation}
    \begin{split}
        f(\infty, y) &> 0 \text{ for } y>0\\
        f(x, \infty) &> 0 \text{ for } x>0.\\
    \end{split}
\end{equation}
There exist polynomials $f(x,y)$ which satisfy (10), but do not satisfy the inequalities in (11).

For example:
\begin{equation*}
    \begin{split}
f(x,y) = x(y-1)^2 + 1\\
f(x,y) \geq 1 \text{ for }x>0,\,y>0\\
f(\infty, y) = (y-1)^2 = 0 \text{ for } y = 1.
    \end{split}
\end{equation*}
For such polynomials, a positive representation is excluded from the outset. If there is any question of this, then one must add the constraints in (12) to the constraints in (11).
We will go somewhat further, and require fulfillment of the following conditions:
\begin{equation}
    \begin{split}
        f(x, y) &> 0 \text{ for } x, y\geq0,\\
        f(\infty, y) &> 0 \text{ for } y\geq0,\\
        f(x, \infty) &> 0 \text{ for } x\geq0,\\
        f(\infty, \infty) &> 0.\\
    \end{split}
\end{equation}

Here, $f(\infty, \infty)$ denotes the term of $f(x,y)$ in which the highest powers of x is multiplied with the highest powers of y. One sees easily that the inequalities in (13) can be homogenized to a single condition.

We now generalize the Laguerrian theorem to the following:

\emph{Theorem: Every polynomial $f(x,y)$ with real coefficients, which satisifes the inequalities in (13), may be written in the positive representation form} 
\begin{equation*}
     f(x,y) = \dfrac{f_2(x,y)}{f_1(x,y)}.
\end{equation*}

The polynomial $f(x,y)$ that fulfills the conditions in (13) has degrees $n_1$ and $n_2$ in $x$ and $y$, respectively.

We order it in $y$:
\begin{equation*}
     f(x,y) = a_0(x) y^{n_2} +  a_1(x) y^{n_2-1} + \ldots a_{n_2}(x).
\end{equation*}
Let $x_0$ be any real nonnegative number. We now examine the roots of the following equation in $y$,
\begin{equation*}
     f(x_0,y) = 0.
\end{equation*}
Because $f(x,\infty) = a_0(x)$, one may write, because of (13), that $a_0(x_0) > 0$. There are exactly $n_2$ roots. We denote them with
\begin{equation}
    y_1(x_0), \cdots, y_{n_2}(x_0).
\end{equation}
Because of constraint (13)a, none are positive real or zero.  Now let $x_0$, starting from its original value, go continuously through all the positive values to $x_0^*$, thus the $n_2$ points of the roots (14) also move steadily in the complex plane according to a known theorem. A point never moves to the real positive axis, or to zero or to infinity. However big $x_0^*$ is, there always exists a region $G^*$ described in the previous section, which no root ever leaves. On the other hand, with sufficiently large $x_0$, the $n$ roots (14) approach the roots of $f(\infty, y) = 0$.

Because of constraint (13)d, all $n_2$ roots are finite, and because of (13)b, none are positive or zero.  Again one can create a region $G^{**}$ akin to $G^*$ that holds all the $n_2$ roots and those roots in (14) that have values $x_0 \geq \bar{x}_0$, where $\bar{x}_0$ is sufficiently large and positive.

We let
\begin{equation*}
    x_0^* = \bar{x}_0
\end{equation*}
and construct a region $G_2$, that contains both $G^*$ and $G^{**}$, which always allows the roots in (14) to remain in this region for any values not less than $x_0$.

Following the method in the previous section, we may form the representation $G_2(y)$, which belongs to region $G_2$.

We note that the number of quadratic factors for the polynomial $f(x_0, y)$ is at most $\left[\dfrac{n_2}{2}\right]$. Thus it follows that for each nonnegative value of $x_0$, the following polynomial in $y$ 
\begin{equation*}
    f(x_0, y) \cdot G_2^{\left[\dfrac{n_2}{2}\right]}(y) = A_0(x)y^N + A_1(x)y^{N-1} + \cdots + A_N(x)
\end{equation*}
has only positive coefficients. The functions
\begin{equation*}
A_i(x), \, i = 0,1,\ldots N
\end{equation*}
whose highest degree in $x$ are $n_1$, also satisfy the inequality
\begin{equation*}
    A_i(x) \geq 0 \text{ for }x \geq 0.
\end{equation*}
All the roots of all equations
\begin{equation*}
   A_i(x)=0, \, i = 0,1,\ldots N
\end{equation*}
therefore lie in the interior of a suitably selected area $G_1$ as in Fig.~7. Let $G_1(x)$ be a similar such representation factor. Then all polynomials
\begin{equation*}
A_i(x) \cdot G_1^{\left[\dfrac{n_1}{2}\right]}(x)
\end{equation*}
have only positive coefficients. Therefore, in the expression
\begin{equation*}
f_2(x,y) = f(x, y) \cdot G_1^{\left[\dfrac{n_1}{2}\right]}(x) G_2^{\left[\dfrac{n_2}{2}\right]}(y),
\end{equation*}
the coefficients of all the elements $x^iy^k$ are positive. We set
\begin{equation*}
f_1(x,y) = G_1^{\left[\dfrac{n_1}{2}\right]}(x) G_2^{\left[\dfrac{n_2}{2}\right]}(y)
\end{equation*}
to obtain the sought positive representation of $f(x,y)$,
\begin{equation*}
f(x,y) = \dfrac{f_2(x,y)}{f_1(x,y)}.
\end{equation*}
This proves the generalization of the Laugerrian theorem.

 \end{document}